\documentclass[10pt,a4paper]{amsart}
\usepackage{pictex}
\title[]{
 Wohlfahrt's Theorem  for the Hecke group  $G_5$ }

\author{\tiny  Cheng Lien Lang and  Mong Lung Lang }

\begin{document}

\baselineskip12pt
\keywords{ Hecke groups, Congruence subgroups, Wohlfahrt's Theorem}
\subjclass[2000]{11F06}

\maketitle
\vspace{-0.3in}
\begin{abstract}
 Let $K$ be a  subgroup of the inhomogeneous Hecke group $G_5$ of geometric
 level $r$. Then $K$ is congruence if and only if $K$ contains
  the principal congruence subgroup $G(2r)$. In the case $r\not\equiv 0$ (mod 4),
    $K$ is congruence if and only if $K$ contains
  the principal congruence subgroup $G(r)$.
  \end{abstract}

\section{Introduction}

\subsection{} Let $H_5$ be the homogeneous  Hecke group generated by $S$ and $T$ given as follows
 and let $H(\pi)
=\{ A \in H_5 \,:\, A \equiv  I \mbox{ (mod $\pi$)}\}$
be the principal congruence subgroup associated with $\pi\in \Bbb Z[\lambda]$.
{\small $$T=
\left (
\begin{array}{rr}
1 & \lambda  \\
0 & 1 \\
\end{array}
\right ),\,
S=
\left (
\begin{array}{rr}
0 & 1 \\
-1 & 0 \\
\end{array}
\right )
\eqno(1.1)$$}
where $\lambda= 2$cos$\pi/5$.
 Let $Z =\left <\pm I\right >$ be the centre of $H_5$. The inhomogeneous Hecke group $G_5$
 and its principal congruence subgroup $G(\pi)$  in which
  a matrix and its negative are identified are defined as
  $G_5 = H_5/Z$ and $G(\pi)= H(\pi)Z/Z$.
A subgroup $K$ of $G_5$ is called a congruence subgroup if $G(\pi) \subseteq K$ for some $\pi$.
 Let $K$ be a congruence subgroup of $G_5$.
 Since $\Bbb Z[\lambda]$ is a principal ideal domain, there exists  a unique
subgroup $G(\pi)$ of $ G_5$
  such that $G(\tau) \subseteq G(\pi)$ whenever $G(\tau) \subseteq K$.
  The ideal $(\pi)$ is called {\em the algebraic level} of $K$.
The  least common multiple $N$ of the cusp widths  of $K$
is called  the {\em geometric level} of $G$
(the width of a cusp $x$ is the
  smallest positive integer   $m$ such that $ \pm T_q^m$    is conjugate
  in $G_5$ to an element of $K$ fixing $x$).
   The main purpose of the present article is to
  prove
 the following theorem  which is an obvious
  generalisation of Wohlfahrt's Theorem for the modular group $PSL(2, \Bbb Z)$.

\smallskip
\noindent {\bf Theorem 4.1.} {\em Let $K\subseteq G_5$ be a congruence subgroup of geometric level $m$
and algebraic level $(\pi)$. Let $n$ be the rational integer below $(\pi)$. Then the following
 holds.
\begin{enumerate}
 \item[(i)] Suppose that $m\not\equiv 0\,\,(mod\,\,4)$. Then $G(m) \subseteq K$ and $n=m$.
  \item[(ii)] Suppose that $4|m$. Then $G(2m)\subseteq K$ and $n$ is either
  $m$ or  $2m$.
  \end{enumerate}
  }
\subsection{} The Plan. Let $m$ be the geometric level of $K$ and $n$  the rational integer
 below the algebraic level $(\pi)$ of $K$. (i) Theorem 4.1
 is equivalent to  $N(G(mn), T^m) = G(m)$ and (ii)  of  Theorem 4.1
 is equivalent to  $G(2m) \subseteq N(G(mn), T^m) $,
  where $N(G(mn), T^m)$ is the smallest normal subgroup of $G_5$ that contains
   $G(mn)$ and $T^m$.
 Since    $N(G(mn), T^m)$ is an
      intermediate subgroup of $G(mn) \subseteq G(m)$, it is clear that
     (i)  Theorem 4.1 can be proved if
       a {\em small} upper bound  $U$ for  $[G(m) : G(mn)]$
       and a {\em large } lower bound $L$ for $[N(H(mn), T^m) : G(mn)]$
        enjoys the fact that $U=L$. This is achieved  in Section 2 (Lemmas 2.1 and 2.3 for the upper bound)
         and Appendices A-C (for the lower bound). As for (ii) of Theorem 4.1, one
          may obtain the fact $G(2m) \subseteq N(G(mn), T^m)$ by studying {\em all}
           the conjugates of $T^m$ modulo $G(mn)$ (Lemmas 3.7 and 3.8).

       \subsection{}
       Geometry of $G_5$.
       An easy study of the group $PSL(2, \Bbb Z[\lambda])$ will
        provide us a good upper bound for $[G(m):G(mn)]$ (Lemma 2.1) except for the case $[G(2\pi) : G(4\pi)]$, where gcd$\,(2,  \pi)=1$.
         This is where the study of the geometry of $G_5$  comes in.
          The actual index of  $[G(2\pi) : G(4\pi)]$ is obtained by studying
          a special polygon (fundamental domain) associated to   $G(2)$ (subsection 2.2).
          As an application of Theorem 4.1, one can prove easily the following proposition.

\smallskip
  \noindent {\bf Proposition 5.1.} {\em Let $K$ be a subgroup of finite index of $G_5$
   with geometric level $m$.
   Then
   \begin{enumerate}
   \item[(i)] Suppose that $m \not\equiv 0\,\,(mod\,\,4)$. Then
   $ K$ is congruence if and only
   if $G(m) \subseteq K$.
   \item[(ii)] Suppose that $4|m$. Then
   $K$ is congruence if and only
   if $G(2m) \subseteq  K$.
   \end{enumerate}}

\subsection{} An Algorithm.
 Wohlfahrt's Theorem has been recognised as one of the indispensable part
 of the study of the congruence subgroups for $SL(2, \Bbb Z)$ and
  $PSL(2,\Bbb Z)$  as
 the geometric level of a subgroup of $PSL(2, \Bbb Z)$ or
 $SL(2,\Bbb Z)$ can be determined
  geometrically (see (v) of subsection 2.2).
 Algorithms for the determination of whether a subgroup of $PSL(2, \Bbb Z)$ is congruence
  can be found in   [LLT2]. With the help of Proposition 5.1,
   the algorithm given in [LLT2] can be generalised
   easily to the Hecke group $G_q$ $(q$ prime). We will investigate subgroups of
    $G_5 $ of index $\le 5$ in Section 5.

\subsection{}The rest of the article is organsied as follows.  In Section 2,
 we give some technical lemmas which will be used to prove Theorem 4.1.
 Section 3 is devoted to
 the study of the normal closure of $G(mn)$ and $T^m$.
 Proof of Wohlfahrt's Theorem for $G_5$ can be found in Section 4.
  Section 5 studies the congruence subgroups of small indices for $G_5$.
   A  proof of Wohlfahrt's Theorem for the modular group
    without Dirichlet' Theorem
    can be found in Appendix D.

\section{Technical Lemmas}
 Let $x,y \in \Bbb Z[\lambda]$.
 Denoted by $(x)$ the ideal generated by $x$.
 We say $x$ is a divisor of $y$
($y$ is a multiple of $x$) if $(y) \subseteq
 (x)$.  For our convenience, for any $x, y\in \Bbb Z[\lambda]$, we will use
  the notation $x|y$ if $x$ is a divisor of $y$.
  Recall that
 \begin{enumerate}
 \item [(i)] $5
 =\lambda^{-2}(2+\lambda)^2
 $ ramifies totally
  in $\Bbb Z[\lambda]$,
  \item[(ii)] 2 and rational primes of the form $10k\pm 3$ are primes in $\Bbb Z[\lambda],$
  \item[(iii)] rational prime of the form $10k\pm 1$ splits into $p = p_1p_2$, where
   $p_i\in \Bbb Z[\lambda]$ are primes.
   \end{enumerate}
 The smallest positive rational integer $m$ in  $(\pi)$ is called the
   rational integer below $(\pi)$.

\subsection{}
For $A$ an ideal of $\Bbb Z[\lambda]$,
we may also define the  principal congruence subgroup $L ( A)$  of $L_5 = SL(2,\Bbb Z[\lambda])$ analogously.
The formula for the index of the principal congruence subgroup $L(A)$ in $SL(2, \Bbb Z[\lambda])$
 is easily calculated as the modular group case (see [Sh]);  it is
$$[L_5 : L( A)] =
%\epsilon
N(A)^3\prod_{P|A}(1-N(P)^{-2}),\eqno(2.1)$$
%\noindent where $ N(A)$ denotes the absolute norm of  $A$ in $\Bbb Z[\lambda]$
% and $\epsilon =1$ if $ A = (2)$, $\epsilon = 1/2$ otherwise.
 For any $u,v \in \Bbb Z[\lambda]$, it is clear that $[G(u) :G(uv)] =\epsilon  [H(u) :H(uv)]$,
 where $\epsilon = 1/2$ or 1 and $\epsilon = 1/2$ if and only if $(2) =(u) \ne (uv)$.
  Applying $(2.1)$, the following is clear.

 \smallskip
\noindent {\bf Lemma 2.1.} {\em Let $n\in \Bbb N$ and let $\pi, \tau  \in \Bbb Z[\lambda]$.
 Suppose that $\pi$ is a prime and that $\tau\notin (\pi)$.
  Then }
  $$[G(\pi \tau ) :G(\pi^{2}\tau)]\le
 \epsilon [H(\pi  ) :H(\pi^{2})]
  \le \epsilon [L(\pi) :L(\pi^{2})] =\epsilon  N(\pi)^3.\eqno(2.2)$$
  {\em where  $\epsilon = 1/2$ or $1$ and $\epsilon = 1/2$ if and only if $(\pi \tau) = (2)$.}

\smallskip
  \noindent {\bf Remark.} $N(2)=4$, $N(2+\lambda)=5.$
 If $p$ is a rational prime, then $N(p) = p^2$. In the case the rational prime below $(\pi)$
  is of the form $p = 10k\pm1 $, $N(\pi)=p$.
\smallskip

\noindent {\bf Lemma 2.2.} {\em
 Let $m$ be a rational integer
and let $p\in \Bbb N$
 be a rational prime divisor of $m$.
 Suppose that $m p \ne  4$.
  Then the following matrices generate an elementary abelian group of
  order $p^6$  modulo $m p$.}
    {\tiny $$\left (\begin{array}{cc}
   1 & m\lambda\\
   0 & 1
 \end{array}\right ),
 \left (\begin{array}{cc}
   1 & 0\\
  - m\lambda & 1
 \end{array}\right ),
 \left  (\begin{array}{cc}
   1-m\lambda^2  & m \lambda^3\\
   -m\lambda& 1+m\lambda^2
 \end{array}\right ),
 \left (\begin{array}{cc}
   1 & m\\
   0 & 1
 \end{array}\right ),
 \left (\begin{array}{cc}
   1 & 0\\
  - m & 1
 \end{array}\right ),
 \left  (\begin{array}{cc}
   1-m \lambda & m\lambda^2 \\
  - m & 1+m\lambda
 \end{array}\right ).$$}
\noindent {\em Proof.}
 Put the above six matrices into the form $I +mU$. One sees easily the following identity.
 $$(I+mU)(I+mV)\equiv I +m(U+V)\equiv
 (I+mU)(I+mV)
 \,\,\, (\mbox{mod} \,\,m p).\eqno(2.3)$$
 Hence
 every non-identity element has order $p$ modulo $mp$ and
 the above matrices  modulo $m p $ generate an abelian group. Note that
  $(2.3)$  makes the unpleasant multiplication of $I +mU$ and  $I+ m V$ into the very easy
  addition of
   $U$ and $V$.
  In order to show the above matrices generate a group of order $p^6$ modulo $m p$, we
   consider the following groups.
{\small $$ M =  \left < X_i =
  \left (
\begin{array}{cc}
1 &  m \lambda ^i   \\
0 & 1 \\
\end{array}
\right ),  Y_i = \left (
\begin{array}{cc}
1 &  0\\
-m\lambda  ^i& 1 \\
\end{array}
\right )
 \right > ,\,\,
N = \left < Z_i =
\left (
\begin{array}{cc}
1 -m\lambda  ^{i+1}     &  m\lambda   ^{i+2}\\
 - m \lambda  ^i& 1 + m \lambda ^{i+1}\\
\end{array}
\right )\right >, $$}
where $0\le i\le 1$.
It is easy to see that
 $M$ and $N$ are abelian groups of order $p^4$ and $p^2$ respectively.
Applying $(2.3)$ and the fact that $N$ is abelian,  elements in $N$ take the following simple form
{\small $$
Z_0^{c_0} Z_1^{c_1} \equiv
\left (
\begin{array}{lr}
1 -m\sum _{i=0}^{1} c_i\lambda^{i+1}    &   m\sum _{i=0}^{1} c_i\lambda^{i+2}                   \\
-m\sum _{i=0}^{1} c_i\lambda  ^{i}  & 1 + m\sum _{i=0}^{1} c_i\lambda^{i+1}  \\
\end{array}
\right ).\eqno(2.4)$$}
Note that we may assume that $0\le c_i\le p-1$.
 Similar to $(2.4)$, elements in $M$ take the form
 {\small  $$  X_0^{a_0} X_1^{a_1}Y_0^{b_0}Y_1^{b_1}\equiv\left (
\begin{array}{cc}
1    &  m\sum _{i=0}^{1} a_i\lambda ^{i}                 \\
-m\sum _{i=0}^{1} b_i\lambda ^{i}  & 1 \\
\end{array}
\right ).\eqno(2.5)$$}
Suppose that
$ Z_0^{c_0} Z_1^{c_1} \equiv  \pm X_0^{a_0} X_1^{a_1}Y_0^{b_0}Y_1^{b_1}$ modulo $mp$
(in $G_5$, a matrix is identified with its negative).
 An easy study of the (2,2)-entries of $(2.4)$ and $(2.5)$
 implies that
 $ 1 + m\sum _{i=0}^{1} c_i\lambda^{i+1}\equiv \pm 1$ (mod $m p)$.
  Since $m p  \ne  4$, one must have $c_0 =c_1=0$. As a consequence,
$M \cap N=\{1\}$. Hence
 $|\left < \Omega \right > | = |M||N| = p^6.$\qed

\subsection{}
In [K], Kulkarni applied a combination of geometric
 and arithmetic methods to show that one can produce
 a set of independent generators in the sense of
 Rademacher for the congruence subgroups of the modular group,
 in fact for all subgroups of finite indices.
 His method can be generalised to all subgroups of
 finite indices of the inhomogeneous Hecke groups $G_q$,
 where $q$ is a prime.
  See [LLT1] for detail
   (Propositions 8-10 and Section 3 of [LLT1]).
    In short,
  for each subgroup $K$ of finite index of $G_q$,
one can associate with $K$
a set of Hecke-Farey symbols (HFS)
$\{-\infty, x_0, x_1, \cdots, x_n, \infty\},$
 a special  polygon (fundamental domain) $\Phi$,
 % and
 %a set of side parings (of this fundamental polygon) which is  a set of
 %independent generators of $S$ (determined by
 % Propositions 8-10 of [LLT1]).
 % To be more precise,
 %every subgroup $G$ of finite index
 %of $G_q1$ admits a generalised Hecke Farey Symbol (gHFS)
 %$$\{ -\infty , x_0, x_1, \cdots , x_n,   \infty \}$$
 and  an additional structure
 on each consecutive pair of $x_i$'s of the three types described
 below :
$$ {x_i} _{_{_{\smile}}} \  \hspace{-.37cm} _{_{_{_{_{_{_{\circ}}}}}}}
  x_{i+1} ,\,\,
{x_i} _{_{_{\smile}}} \  \hspace{-.37cm} _{_{_{_{_{_{_{\bullet}}}}}}}
  x_{i+1} ,\,\,
{x_i} _{_{_{\smile}} }\  \hspace{-.37cm} _{_{_{_{_{_{_{a}}}}}}}  x_{i+1}.
 $$
where $a$ is a nature  number. Each nature number $a$ occurs
 exactly twice or not at all.
Similar to the modular group,
 the actual values of the $a$'s is unimportant: it is
 the pairing induced on the  pairs that matters.
 \begin{enumerate}
 \item[(i)] The side pairing  $\circ$ is an elliptic element of order 2 that pairs
  the even line $(a/b, c/d)$ with itself. The trace of such an element is 0.
    \item[(ii)]
  The side pairing  $\bullet$ is an  elliptic element of order $q$ that pairs
  the odd line $(a/b, c/d)$. The absolute value of the trace of such an element is
   $\lambda_q$.
      \item[(iii)] The two sides $(a/b, c/d)$ and $(u/v, x/y)$ with the label $a$
   are paired together by an  element of infinite order.
  \item[(iv)]
The special polygon associated with the HFS is a fundamental
 domain of $K$ and the side pairings
 $ I_K = \{ g_1,  g_2 \cdots ,  g_m \}$
associated with
 the  HFS
 is a set of independent generators of $K$
 (Theorem 7, Propositions 8-10 of [LLT1]).
% \item[(v)] The special polygon is a union of ideal $q$-gons and special
%  triangles. Each $q$-gon is a union of $q$ special triangles.
%  \item[(vi)]
%The number $d$ of  special triangles
% (a special triangle is  a fundamental domain of $G_q$) of the special polygon is  the index
% of the
%   subgroup.
%   \item[(vi)] The set of independent generators consists of
%  $r$ matrices of infinite order, where $r$ is the number of the nature number $a$'s in the
%   Hecke-Farey smybols.
%   \item[(vii)] The subgroup has
%   $v_2$
%   (the number of the circles $\circ$ in HFS) inequivalent classes of elliptic
%    elements of order 2. Each class has exactly one representative in $I$.
%     \item[(viii)] The subgroup has
%   $v_q$
 %  (the number of the bullets $\bullet$ in HFS) inequivalent classes of elliptic
%    elements of order $q$. Each class has exactly one representative in $I$.
  \item[(v)] The width of a cusp $x$,
      denoted by $w(x)$,  is the number of even lines
      in $\Phi$ that comes into $x$. Algebraically, it is the smallest positive integer
       $m$ such that $\pm T_q^m$ is conjugate in $G_q$ to an element of $K$ fixing $x$
       (keep in mind that a matrix is identified with its negative in $G_q$).                            \end{enumerate}

 \smallskip
 \noindent {\bf Lemma 2.3.} {\em
 $G(2) / G(4)$ is an elementary abelian group of order $2^4$. Suppose that gcd$\,(2,\pi)=1$ and
  $(\pi)\ne(1)$.
  Then
 $G(2\pi) / G(4\pi)$ is an elementary abelian group of order $2^5$.
 }
 \smallskip

 \noindent {\em Proof.}
 A set of independent generators of $G(2)$ is given as follows
  (see [LLT1] for detail).
{\small $$\Omega_2 =\left \{\left (\begin{array}{cc}
   1 & 2\lambda\\
   0 & 1
 \end{array}\right ),
 \left (\begin{array}{cc}
   1 & 0\\
   2\lambda & 1
 \end{array}\right ),
 \left (\begin{array}{cc}
   1+2\lambda & 2+ 2\lambda\\
   2\lambda & 1 +2\lambda
 \end{array}\right ),
 \left (\begin{array}{cc}
   1+2\lambda & 2\lambda\\
   2 +2\lambda & 1+2\lambda
 \end{array}\right )\right \}.\eqno(2.6)$$}
 Denoted by $a_i$ the members in $\Omega _2$. It is easy to see that
  $a_i^2 \equiv 1$ and $a_ia_j\equiv a_ja_i$ modulo 4.
  Similar to Lemma 2.2, one can show
  that $[G(2) : G(4)] = 2^4$. By Lemma 2.1,
   $[G(2\pi) : G(4\pi)] \le
   [H(2) :H(4) ] =2^5$.
   Let $m$ be the smallest rational integer in $(\pi)$.
   By our results in Appendix B, $G(2\pi)$ contains the following five matrices
    (modulo $4\pi$).
  {\tiny $$  I+ 2m\left (\begin{array}{cc}
   0 & \lambda\\
   0 & 0
 \end{array}\right ), I+2m
 \left (\begin{array}{cc}
   0 & 0\\
  -\lambda & 0
 \end{array}\right ), I+2m
 \left  (\begin{array}{cc}
   -\lambda & 0\\
   -1& \lambda
 \end{array}\right ),  I+2m
 \left (\begin{array}{cc}
   0 & 1\\
   -1 & 0
 \end{array}\right ),  I+2m
  \left (\begin{array}{cc}
   -1  & 1 \\
   -1 & 1
 \end{array}\right ).
 $$}
    Since these matrices  generate an elementary
    abelian group of order $2^5$ (see Appendix B), one has
    $[G(2\pi) : G(4\pi)] =2^5$.
  \qed

\subsection{} The homogeneous Hecke group $H_5$.  In $G_5$, $G(ab) \ne G(a)\cap G(b)$.
 This  creates some obstruction for our study of the normal closure (see Lemma 3.2)
  and we have to turn our study to the homogeneous group $H_5$ as follows.

 \smallskip
 \noindent {\bf Lemma 2.4.} {\em
  Let $K$ be a normal subgroup of $H_5$. Suppose
  that $T\in K$. Then $K= H_5$.
  Let $a, b\in \Bbb N$.
  Suppose that gcd$\, (a, b)=1$. Then $H(a)H(b) = H_5$ and
   $H_5/H(ab)\cong H_5/H(a)\times H_5/H(b)$.}

  \smallskip
  \noindent
   {\em Proof.}
   We note first that we are working on $H_5$, where $H(a)\cap H(b) = H(ab)$
    and the order of $S$ is 4 $(S$ has order 2 in $G_5)$.
     Every element $x$  in $H_5$ can be written as a word $w(S, T)$ in $S$ and $T$.
    Since $K$ is normal, $T \in K$, and $S$ is of order 4, we have $xK
     = K$, $SK$, $S^2K$ or $S^3K$.
      Hence $[G_5 : K] = 1$, 2 or 4. Note that $ST$ is of order $5$.
       Hence $ST = (ST)^{16} \in K$.  Since $T\in K$,
     we have  $S\in K$. As a consequence, $K=H_5$.
       We shall now study the second part of the lemma.
       Since gcd$\,(a, b)=1$, $T^a\in H(a)$, and $T^b\in H(b)$, one has $T\in H(a)H(b)$.
       This implies that $H_5 = H(a)H(b)$.
        Note that $H(a)\cap H(b) = H(ab)$.
      As a consequence,
       $H_5/H(ab)\cong H_5/H(a)\times H_5/H(b)$.\qed

\smallskip

Lemma 2.4 can be generalised to all Hecke group $G_q$, where $q$ is a prime, as follows.

 \smallskip
 \noindent {\bf Lemma 2.5.} {\em Let $q$ be a prime and let $K$ be a normal subgroup of $G_q$. Suppose
  that $T\in K$. Then $K= G_q$.
  In particular, for any $a, b\in \Bbb N$, if gcd$\, (a, b)=1$, then $G(a)G(b) = G_q$.
  }

\section {The normal closure   $N(G(\pi), T^m)$    }

 \noindent {\bf Definition 3.1.} Let $m\in \Bbb N$ and $ \pi\in \Bbb Z[\lambda].$
  Denoted  by $N(G(\pi), T^m)$ the smallest normal subgroup
  (normal closure) of $G_5$ that contains $G(\pi)$
   and $T^m$.

\subsection{} The main purpose of this subsection is to study the group $N(G(mn), T^m)$
 where $m$ or $n$ is odd, or gcd$\,(m,n)=1$. We will show that
 $N(G(mn), T^m) = G(m)$ under such assumption.

 \smallskip
\noindent {\bf Lemma 3.2.}   {\em Let  $a, b\in \Bbb N$. Suppose that
 gcd$\,(a,b)=1$. Then  $N(G(ab)), T^b)=G(b)$. }

 \smallskip
 \noindent {\em Proof.}  To prove our lemma, we shall work on $H_5$, where
  $H(a)\cap H(b) = H(ab)$ and $H_5\cong H(b)/H(ab)\times H(a)/H(ab)$ (Lemma 2.4).
 Set $X = N(H(ab), T^b)$. Then
 $H(ab) \subseteq X \subseteq H(b)$.
  Since $T^b\in X, T^a\in H(a)$ and gcd$\,(a,b)=1,$
  one has $T\in XH(a)$. By Lemma 2.4, $XH(a) = H_5$. Hence
  $$
  {\small X/H(ab) \times H(a)/H(ab) \cong XH(a)/H(ab)=H_5/H(ab) \cong H(b)/H(ab)\times H(a)/H(ab).}\eqno(3.1)$$
   Note that $X \subseteq H(b)$.
  It is now clear that (3.1) implies that $X = H(b)$. Equivalently,
   $ H(b)=  N(H(ab), T^b)$.
    As a consequence, one has $G(b) = N(G(ab), T^b)$. \qed
  \smallskip

  \noindent {\bf Remark } In Lemma 3.2,  one has to work on $H_5$ rather than $G_5$
   as $H(a)\cap H(b) = H(ab)$ but $G(a)\cap G(b) \ne G(ab)$.

   \smallskip
\noindent {\bf Lemma 3.3.} {\em Suppose that $n\in \Bbb N$ is odd.
  Them $N(G(mn), T^m) = G(m)$.}

 \smallskip
\noindent {\em Proof.}
We shall first assume that
$n=p$ is an odd prime.
    Applying $(A8)$ of  Appendix A, $  N(G(mp), T^m)$ contains the following
      six matrices (modulo $mp$).
   {\tiny $$  \left (\begin{array}{cc}
   1 & m\lambda\\
   0 & 1
 \end{array}\right ),
 \left (\begin{array}{cc}
   1 & 0\\
  - m\lambda & 1
 \end{array}\right ),
 \left  (\begin{array}{cc}
   1-m\lambda^2  & m\lambda^3\\
   -m\lambda& 1+m\lambda^2
 \end{array}\right ),
 \left (\begin{array}{cc}
   1 & m\\
   0 & 1
 \end{array}\right ),
 \left (\begin{array}{cc}
   1 & 0\\
  - m & 1
 \end{array}\right ),
 \left  (\begin{array}{cc}
   1-m \lambda & m\lambda^2 \\
  - m & 1+m\lambda
 \end{array}\right ).$$}
\hspace{-.2cm} By Lemma 2.2, $|N(G(mp), T^m)/G(mp)|
    \ge p^6$.
    Since $N(G(mp), T^m)/G(mp)\subseteq G(m)/G(mp)$
     and $|G(m)/G(mp)| \le p^6$ (Lemma 2.1),
     we have $
      N(G(mp), T^m)= G(m)$.
We shall now study the general case.
It is clear that $N(G(nm), T^m) \subseteq  G(m)$.
 Let $n_0$ be the smallest positive integer such that $G(n_0 m)\subseteq G(m)$.
  Suppose that $n_0 >1$. Let $p$ be a prime divisor of $n_0$ Then
  $$G(n_0m/p) =
   N(G(p(n_0m/p)), T^{mn_0/p}) \subseteq N(G(nm), T^m) \subseteq G(m).\eqno(3.2)$$
   This contradicts the minimality of $n_0$. It follows that $n_0=1$ and that $N(G(nm), T^m) = G(m)$.
   This completes the proof of the lemma.
    \qed

 \smallskip
    \noindent {\bf Lemma   3.4.} {\em
    Let $m, n\in \Bbb N$.
   Suppose that $m$ is odd.
   Then  $N(G(nm), T^m)=G(m)$.}

  \smallskip
\noindent {\em Proof.}
 Decompose $n$ into $n = m' m_0$, where $m_0=$ gcd$\,(m,n)$.
  It follows that gcd$\,( m', m_0m)=1$.
  It is clear that
 $$N(G(m'm_0 m), T^{m_0m}) \subseteq N(G(nm), T^m).\eqno(3.3)$$
 Since gcd$\,( m', m_0m)=1$, we may apply Lemma 3.2 and conclude that the group to
  the right
  hand side of (3.3) is
 $$N(G(m'm_0 m), T^{m_0m})=G(m_0m).\eqno(3.4)$$
  Hence $G(m_0m) \subseteq N(G(nm), T^m)$. As a consequence,
  $N( G(m_0m), T^m ) \subseteq N(G(nm), T^m).$
 Since $m$ is odd and   $m_0=$ gcd$\,(m,n)$, $m_0$ is odd as well.
  Applying Lemma 3.3, one has
  $$G(m) = N( G(m_0m), T^m )
  \subseteq   N(G(nm), T^m)\subseteq G(m)
  .\eqno(3.5)$$
 Hence $N(G(nm), T^m) = G(m)$.
 This completes the proof of the lemma.
 \qed

\subsection{}
The main purpose of this subsection is to study the group $N(G(mn), T^m)$ when $n$ or $m$ is even.
It turns out that the group $N(G(mn), T^m)$  does not behave very well  as their counterparts
in subsection 3.1.

  \smallskip
    \noindent {\bf Lemma 3.5.} {\em Let $m\in \Bbb N$. Suppose that gcd$\,(m, 2) =1$. Then
     $N(G(4m), T^{2m} )=G(2m)$.}

     \smallskip
     \noindent {\em Proof.}
      By our results in Appendix B,   $N(G(4m), T^{2m} ) \subseteq G(2m) $ possesses the following
      five  matrices (modulo $4m$).
  {\tiny $$  I+ 2m\left (\begin{array}{cc}
   0 & \lambda\\
   0 & 0
 \end{array}\right ), I+2m
 \left (\begin{array}{cc}
   0 & 0\\
  -\lambda & 0
 \end{array}\right ), I+2m
 \left  (\begin{array}{cc}
   -\lambda & 0\\
   -1& \lambda
 \end{array}\right ),  I+2m
 \left (\begin{array}{cc}
   0 & 1\\
   -1 & 0
 \end{array}\right ),  I+2m
  \left (\begin{array}{cc}
   -1  & 1 \\
   -1 & 1
 \end{array}\right ).
 $$}
  The above matrices modulo $4m$  generate an elementary abelian group of
  order $2^5$ if $m\ge 2$, an elementary abelian group of order
   $2^4$ if $m =1$ (see Appendix B). Note that the orders meet
   the upper bound of $[G(2m) : G(4m)]$ is both cases (Lemma 2.3).
    As a consequence, one has
     $N(G(4m), T^{2m}) = G(2m)$.\qed

   \smallskip
   \noindent {\bf Lemma 3.6.} {\em  Let $m\in \Bbb N$. Suppose that gcd$\,(m, 2) =1$. Then
     $N(G(8m), T^{2m} )=G(2m)$.}

   \smallskip
   \noindent {\em Proof.} It is clear
     that $N(G(8m),  T^{4m} ) \subseteq N(G(4m), T^{2m})$.
      Apply our results in Appendix B, the following are members of $N(G(8m),  T^{4m} )
      \subseteq G(4m) \cap N(G(8m), T^{2m})$ modulo $8m$.
  {\tiny $$ I+ 4m\left (\begin{array}{cc}
   0 & \lambda\\
   0 & 0
 \end{array}\right ) , I+4m
 \left (\begin{array}{cc}
   0 & 0\\
  -\lambda & 0
 \end{array}\right ), I+4m
 \left  (\begin{array}{cc}
   -\lambda & 0\\
   -1& \lambda
 \end{array}\right ), I+4m
 \left (\begin{array}{rc}
   0 & 1\\
   -1 & 0
 \end{array}\right ),  I+4m
  \left (\begin{array}{cc}
   -1  & 1 \\
   -1 & 1
 \end{array}\right ).
 $$}
 Denote these five matrices by $X_1, X_2, X_3, X_4$ and $X_5$.
Then $\{ X_1, X_2, X_3, X_4, X_5\}$ generates an elementary abelian group  $E$ of order $2^5$ modulo $8m$ (see Appendix B).
Note that $E \subseteq G(4m)/G(8m)$ and that
 $$
 {\tiny  I+4m
  \left (\begin{array}{rc}
   0   & 0 \\
   1 &0
 \end{array}\right )} \mbox{ is not an element of }  E  \subseteq G(4)/G(8m).\eqno(3.6)
 $$
 We now consider the element $ (T^2ST^{-2}S^{-1})^2$.
 Since $4m^2 \equiv 4m$ (mod $8m$), $ \lambda ^2 = \lambda +1$,  direct calculation shows that
  (modulo $8m$)
 $$ {\tiny  X=   (T^2ST^{-2}S^{-1})^2
 \equiv I+ 4m
  \left (\begin{array}{cc}
   \lambda +1  & \lambda \\
   \lambda  &  \lambda+1
 \end{array}\right )} \in N(G(8m), T^{2m}) \cap G(4m).\eqno(3.7)
 $$
 Applying the identity $(I+4mU)(I+4mV) \equiv I+ 4m(U+V)$ modulo $8m$, the following is clear.
    $$X_1X_2X_3X_4X_5X \equiv
  {\tiny  I+4m
  \left (\begin{array}{rc}
   0   & 0 \\
   1 &0
 \end{array}\right )}\eqno(3.8)$$
 By (3.6), $X_1X_2X_3X_4X_5X \notin E$. Hence $\left <X_1, X_2, X_3, X_4, X_5, X\right >$
  modulo $8m$
  is  an elementary abelian group of order $2^6$. Note that
   $$ F =\left <X_1, X_2, X_3, X_4, X_5, X\right > \subseteq G(4m)/G(8m) \cap N(G(8m),T^{2m})/G(8m).
   \eqno(3.9)$$
    Since $|F|=2^6$ and $[G(4m) : G(8m)]\le 2^6$ (Lemma 2.1), we conclude that
     $F= G(4m)/G(8m)$ and that $G(4m)\subseteq N(G(8m), T^{2m})$.
      As a consequence, $N(G(4m), T^{2m}) \subseteq N(G(8m), T^{2m}) \subseteq G(2m)$.
       By Lemma 3.5, $N(G(4m), T^{2m}) =G(2m)$. Hence
    $ N(G(8m), T^{2m}) = G(2m)$.
    This completes the proof of the lemma.\qed

 \smallskip
\noindent {\bf Lemma 3.7.} {\em Let $m\in \Bbb N$.
 Suppose that $m$ is a multiple of $4$. Then $ G(2m) \subseteq N(G(4m), T^m)$.
   }

 \smallskip
 \noindent
 {\em Proof.}
   By $(C4)$ of  Appendix C, $  N(G(4m), T^m)$ contains the following
      six matrices (modulo $4m$).
        {\tiny $$ \left (\begin{array}{cc}
   1 & 2m\lambda\\
   0 & 1
 \end{array}\right ),
 \left (\begin{array}{cc}
   1 & 0\\
  - 2m\lambda & 1
 \end{array}\right ),
 \left  (\begin{array}{cc}
   1-2m\lambda^2  & 2m\lambda^3\\
   -2m\lambda& 1+2m\lambda^2
 \end{array}\right ),
 \left (\begin{array}{cc}
   1 & 2m\\
   0 & 1
 \end{array}\right ),
 \left (\begin{array}{cc}
   1 & 0\\
  -2 m & 1
 \end{array}\right ),
 \left  (\begin{array}{cc}
   1-2m \lambda & 2m\lambda^2 \\
  - 2m & 1+2m\lambda
 \end{array}\right ).$$}
\hspace{-.2cm} Applying  Lemma 2.2,
the above matrices generate an  abelian group of $G(2m)/G(4m)$ of order $2^6$.
 Since  $|G(2m)/G(4m)| \le 2^6$ (Lemma 2.1), we have
    $G(2m) \subseteq N(G(4m), T^m)$. \qed

\smallskip
\noindent {\bf Lemma 3.8.} {\em Let $m \in \Bbb N$. Suppose that $2|m$. Then
 $[N(G(4m), T^{2m}): G(4m)]=2^5$.
 In particular, $N(G(4m), T^{2m})\ne G(2m)$.
 }

 \smallskip
 \noindent {\em Proof.} We note first that although
 members in $(C4)$ are not generated by the conjugates of $T^{2m}$ but they are
  members of  $G(2m)$. As a consequence, one has
  $[G(2m) : G(4m)] \ge  2^6$. This implies that   $[G(2m) : G(4m)] =  2^6$ (Lemma 2.1).
  Let $S_1 = \{T^{2m}\}$ and define  $S_{n+1}$
  inductively to be $S_{n+1}  = S_n\cup \{ Sx_iS^{-1}, Tx_iT^{-1} \,:\, x_i \in S_n\}$.
   It is clear that $S_i \subseteq G(2m)$ for all $i$.
Each $S_i$ generates an elementary abelian group modulo $G(4m).$
 Since $[G(2m) : G(4m) ] = 2^6$ is finite, there exists some $r$ such that
  $\left < S_r\right > =\left < S_{r+1}\right >$ modulo $G(4m)$.
  It follows  that $N(G(4m), T^{2m}) = \left < S_r\right >G(4m)/G(4m)$.
   Direct calculation shows that $\left < S_5 \right > =\left < S_{6}\right >
   = \left < S_{7}\right > = \cdots
   $ is given as follows.
    {\small $$\left < \left (\begin{array}{cc}
   1 & 2m\lambda\\
   0 & 1
 \end{array}\right ),
 \left (\begin{array}{cc}
   1 & 0\\
  2m\lambda & 1
 \end{array}\right ),
 I+2m
 \left  (\begin{array}{cc}
   1  &0\\
   0& 1
 \end{array}\right ), I+2m
 \left (\begin{array}{cc}
   \lambda+1 & 1\\
   0 & \lambda+1
 \end{array}\right ), I+2m
 \left (\begin{array}{cc}
   0 & 1\\
   1 & 0
 \end{array}\right )\right >.
 $$}
One can show by direct calculation that the above five matrices modulo $4 m$
 generate an elementary abelian group of order $2^5$.
  Hence $[N(H(4m), T^{2m}) : G(4m)] = 2^5$ and
   $N(G(4m), T^{2m}) \ne G(2m)$.\qed

 \smallskip
\noindent {\bf Lemma 3.8.} {\em Let $m\in \Bbb N$.
 Suppose that $m$ is a multiple of $4$. Then $ G(2m) \subseteq N(G(4m), T^m)$.
   }

 \smallskip
 \noindent
 {\em Proof.}
   By $(C4)$ of  Appendix C, $  N(G(4m), T^m)$ contains the following
      six matrices (modulo $4m$).
        {\tiny $$ \left (\begin{array}{cc}
   1 & 2m\lambda\\
   0 & 1
 \end{array}\right ),
 \left (\begin{array}{cc}
   1 & 0\\
  - 2m\lambda & 1
 \end{array}\right ),
 \left  (\begin{array}{cc}
   1-2m\lambda^2  & 2m\lambda^3\\
   -2m\lambda& 1+2m\lambda^2
 \end{array}\right ),
 \left (\begin{array}{cc}
   1 & 2m\\
   0 & 1
 \end{array}\right ),
 \left (\begin{array}{cc}
   1 & 0\\
  -2 m & 1
 \end{array}\right ),
 \left  (\begin{array}{cc}
   1-2m \lambda & 2m\lambda^2 \\
  - 2m & 1+2m\lambda
 \end{array}\right ).$$}
\hspace{-.2cm} Applying  Lemma 2.2,
the above matrices generate an  abelian group of $G(2m)/G(4m)$ of order $2^6$.
 Since  $|G(2m)/G(4m)| \le 2^6$ (Lemma 2.1), we have
    $G(2m) \subseteq N(G(4m), T^m)$. \qed

\smallskip
\noindent {\bf Discussion 3.9.} Lemma 3.7 cannot be improved as
  $N(G(2m), T^m)
 \ne G(m)$ for $4|m$ (Lemma 3.8).

\section{Wohlfahrt's Theorem for $G_5$}

\smallskip
\noindent {\bf Theorem 4.1.} {\em Let $K\subseteq G_5$ be a congruence subgroup of geometric level $r$
and algebraic level $(\pi)$. Let $n$ be the rational integer below $(\pi)$. Then the following
 holds.
\begin{enumerate}
 \item[(i)] Suppose that $r$ is odd . Then $G(r) \subseteq K$ and $n=r$.
 \item[(ii)] Suppose that $r=4k+2$. Then $G(r)\subseteq K$ and $n=r$.
   \item[(iii)] Suppose that $4|r$. Then $G(2r)\subseteq K$ and $n$ is either
  $r$ or  $2r$.
  \end{enumerate}
  }

\smallskip
\noindent {\em Proof.} We shall first study the case that $r$ is odd.
 Since the geometric level is
 $r$, $\pm xT^rx^{-1} \in K$
for all $x \in G_5$.
 Note that $G(n) \subseteq G(\pi) \subseteq K$.
Hence $K$ contains $N(G( nr), T^r)$.
By   Lemma 3.4, $G(r) \subseteq K$. Since the algebraic level is $(\pi)$,
 $G(r) \subseteq G(n)$. Hence $n|r$.  Conversely, since $G(n) \subseteq K$, the
  geometric level of $K$ is a divisor of $n$. Hence $r|n$. It follows that $n=r$.

\smallskip
\noindent
 (ii)
 Since the geometric level is $r$, $\pm xT^rx^{-1} \in K$ for all $x\in G_5$.
This implies that $N(G(nr), T^{r}) \subseteq K$. Set $nr = (2k+1)2^s r$.
 By Lemma 3.3, one has
$$G(2^sr) = N(G((2k+1)2^sr), T^{2^sr}) \subseteq N(G(nr), T^{r})\subseteq K.\eqno(4.1)$$
 As a consequence,
 $$N(G(2^sr), T^{r}) \subseteq N(G(nr), T^{r})
 \subseteq K.\eqno(4.2)$$
 We now study the group
$N(G(2^sr), T^{r}) $ in (4.2).
It is clear that if $s\le 2$. Then one has the following.
$$G(r) = N(G(2^sr), T^{r})\subseteq K \,\,\,\mbox{ (Lemmas 3.5,  3.6)}.\eqno(4.3)$$
We shall now assume that $s \ge 3$. Let $k$ be the smallest positive integer such
 that $G(2^{k} r)\subseteq K$. Suppose that $2^k\ge 8$. By Lemma 3.7,
 $$G(2^{k-1}r) = N(G(4(2^{k-2}r), T^{2^{k-2}r}) \subseteq N(G(2^kr), T^{r})\subseteq K.\eqno(4.4)$$
This contradicts the minimality of $k$. Hence $s\le 2 $ and $G(r)\subseteq K$.
 This completes the proof of (ii).
   (iii) can be proved similarly.
   \qed

\smallskip
\noindent {\bf Discussion 4.2.}
 Unlike the modular group case, the algebraic level $(\pi)$ of $K$ is not uniquely
    determined by the geometric level $r$. For instance,
    while the algebraic level of $G(5)   $ and $G(2+\lambda)  $ are $(5)$ and $(2+\lambda)$
     respectively,
     both $G(5)$ and $G(2+\lambda)$
     are congruence of geometric level $5$.

\smallskip
\noindent {\bf Discussion 4.3.}  (iii) of Theorem 4.1
cannot be improved as $K= N(G(8), T^4)$ has geometric level 4 and algebraic 8 (see Lemma 3.8).

\section{Congruence Subgroup Problem}
Let $K$ be a subgroup of $G_5$ of finite index and let $\Phi$ be a special polygon of $K$
 (see subsection 3.2). Applying (v) of subsection 3.2, the geometric level $r$ of $K$
  can be determined geometrically.
   By Theorem 4.1,
  we have the following.

  \smallskip
  \noindent {\bf Proposition 5.1.} {\em Let $K$ be a subgroup of finite index of $G_5$
   with geometric level $m$.
   Then
   \begin{enumerate}
   \item[(i)] Suppose that $m \not\equiv 0\,\,(mod\,\,4)$. Then
   $ K$ is congruence if and only
   if $G(m) \subseteq K$.
    \item[(ii)] Suppose that $4|m$. Then
   $K$ is congruence if and only
   if $G(2m) \subseteq  K$.
   \end{enumerate}}

\smallskip
An algorithm for the determination of whether a subgroup of $PSL(2, \Bbb Z)$ is congruence
  can be found in   [LLT2].
  With the help of Proposition 5.1, the algorithm given in [LLT2] can be generalised
   easily to the Hecke group $G_q$ $(q$ prime).
    An easy observation of the special polygons implies that
    \begin{enumerate}
   \item[(i)] There is a unique subgroup of index 2 given as in subsection 5.1.
    \item[(ii)] $G_5$ has subgroups of
     all possible indices except for 3 and 4.
     \item[(iii)] A special polygon of a subgroup of index 5
      consists of one special 5-gon. Its side pairings must have an element of order 2.
       There are altogether 26 subgroups of index 5 (see [LLT3]).
      \item[(iv)] The only normal
        subgroup of index 5 is $K = \left < x^5 \,:\, x\in G_5\right >$, a set
         of independent generators of $K$ consists of 5 elements of order 2 (see subsection 5.7).
         \end{enumerate}
     Note that since 5 is a prime, subgroups
     of index 5 are either normal or self-normalised. As a consequence, each non-normal
     subgroup of index 5 has 5 conjugates.
    Our study of subgroups of index $\le 5$ is recorded as follows.

\subsection{} $G_5$ has a unique subgroup of index 2 with the following
 Hecke Farey symbol. It is congruence of geometric level 2. The algebraic level is also 2.
$${\small  -\infty _{_{_{\smile}}} \  \hspace{-.37cm} _{_{_{_{_{_{_{\bullet}}}}}}}
 0
  _{_{_{\smile}}} \  \hspace{-.37cm}
%   _{_{_{_{_{_{_{\circ}}}}}}}
 % 1 _{_{_{\smile}} }\  \hspace{-.37cm}
  _{_{_{_{_{_{_{\bullet}}}}}}}  \infty  }\eqno(5.1)
 $$

 \subsection{} Index 5,  $v_2=1$ with  geometric level 2.
  $G_0(2)$ is one of them with the following Hecke Farey symbols.
   $G_0(2)$ is self normalised. Consequently,
    $G_0(2)$ has 5 conjugates. $ \Bbb Z_2 \cong G_0(2)/G(2) \subseteq G_5/G(2) \cong D_{10}$.
  $${\small
 -\infty _{_{_{\smile}}} \  \hspace{-.37cm} _{_{_{_{_{_{_{_{1}}}}}}}}
 0 \, _{_{_{\smile}}} \  \hspace{-.37cm} _{_{_{_{_{_{_{_{2}}}}}}}}
\, 1/\lambda _{_{_{\smile}} }\   \hspace{-.37cm} _{_{_{_{_{_{_{_{\circ}}}}}}}}
  \lambda/\lambda _{_{_{\smile}} }\  \hspace{-.37cm} _{_{_{_{_{_{_{_{2}}}}}}}}
  \lambda _{_{_{\smile}} }\  \hspace{-.37cm} _{_{_{_{_{_{_{_{1}}}}}}}}  \infty} \eqno(5.2)
 $$
 \subsection{}Index 5,  $v_2=1$ with  geometric level 3.
 The following subgroup  $K$ is one of them.
$K$ has 5 conjugates.
   Since $ G_5/G(3) \cong A_5$ has 5 subgroups of index 5, $G_5$ has  5 congruence
    subgroups of geometric level 3, index 5. It follows that $K$ and its conjugates are
     congruence as they are the only subgroups (of geometric level 3) among the 26 subgroups
      of index 5. The algebraic level of these groups is also 3.
    $${\small
 -\infty _{_{_{\smile}}} \  \hspace{-.37cm} _{_{_{_{_{_{_{_{1}}}}}}}}
 0 \, _{_{_{\smile}}} \  \hspace{-.37cm} _{_{_{_{_{_{_{_{1}}}}}}}}
\, 1/\lambda _{_{_{\smile}} }\   \hspace{-.37cm} _{_{_{_{_{_{_{_{\circ}}}}}}}}
  \lambda/\lambda _{_{_{\smile}} }\  \hspace{-.37cm} _{_{_{_{_{_{_{_{2}}}}}}}}
  \lambda _{_{_{\smile}} }\  \hspace{-.37cm} _{_{_{_{_{_{_{_{2}}}}}}}}  \infty }\eqno(5.3)
 $$

  \subsection{}Index 5,  $v_2=1$ with  geometric level 5. The following subgroup
  $K$ is one of them. $K$ has 5 conjugates.
 Since $ G_5/G(2+\lambda) \cong A_5$ has 5 subgroups of index 5, $G_5$ has  5 congruence
    subgroups of geometric level  5, algebraic level $2+\lambda$, index 5. It follows that $K$ and its conjugates are
     congruence as they are the only non-normal subgroups (of geometric level 5) among the 26 subgroups
      of index 5. The algebraic level of these groups is $2+\lambda$.
  $${\small
 -\infty _{_{_{\smile}}} \  \hspace{-.37cm} _{_{_{_{_{_{_{_{1}}}}}}}}
 0 \, _{_{_{\smile}}} \  \hspace{-.37cm} _{_{_{_{_{_{_{_{2}}}}}}}}
\, 1/\lambda _{_{_{\smile}} }\   \hspace{-.37cm} _{_{_{_{_{_{_{_{\circ}}}}}}}}
  \lambda/\lambda _{_{_{\smile}} }\  \hspace{-.37cm} _{_{_{_{_{_{_{_{1}}}}}}}}
  \lambda _{_{_{\smile}} }\  \hspace{-.37cm} _{_{_{_{_{_{_{_{2}}}}}}}}  \infty }\eqno(5.4)
 $$

\subsection{}Index 5,
$v_2 = 3$ with geometric level 4.
The following subgroup $K$ is one of them. $K$ has 5 conjugates.
 Suppose that $K$ is congruence. By Proposition 5.1, $G(8) \subseteq K$.
 Note that  $[G_5 : G(8) ] = 2^{11}5$, $[G_5 : K]= 5$.
 Hence $K/G(8)$ is a Sylow 2-subgroup of $G_5/G(8)$.
  Note that $G_0(2)/G(8)$ is also a Sylow 2-subgroup of $G_5/G(8)$.
  Hence $K$ and $G_0(2)$ are conjugate to each other. This is a contradiction
   as these two groups have different geometric invariants. Hence $K$ is not
   congruence.
  $${\small
 -\infty _{_{_{\smile}}} \  \hspace{-.37cm} _{_{_{_{_{_{_{_{ 1  }}}}}}}}
 0 \, _{_{_{\smile}}} \  \hspace{-.37cm} _{_{_{_{_{_{_{_{\circ }}}}}}}}
\, 1/\lambda _{_{_{\smile}} }\   \hspace{-.37cm} _{_{_{_{_{_{_{_{\circ }}}}}}}}
  \lambda/\lambda _{_{_{\smile}} }\  \hspace{-.37cm} _{_{_{_{_{_{_{_{\circ }}}}}}}}
  \lambda _{_{_{\smile}} }\  \hspace{-.37cm} _{_{_{_{_{_{_{_{1 }}}}}}}}  \infty }\eqno(5.5)
 $$

\subsection{}Index 5,
$v_2 = 3$ with geometric level 6. $H_5/H(6) \cong H(2)/H(6)\times H(3)/H(6)
 \cong SL(2, 5) \times D_{10}$. Let $A$ be a subgroup of $H(2)/H(6)$ of order $24$ and let
  $B$ be a subgroup of $H(3)/H(6)$ of order 2. Set $V = \left < A\times B, (a,b)\right >$
   where $a$ is an element of $H(2)/H(6)$ of order 5 and
 $b$ is an element of $H(3)/H(6)$ of order 5. It is clear that $H(2)/H(6)$ and $H(3)/H(6)$
  are not subgroups of $V$ and that $V$ is of index 5 in $H_5/H(6)$.
  Set $K = VZ/Z$ ($Z$ is the centre of $H_5$). Then $K$ is a subgroup of $G_5/G(6)$ of index 5.
  Note that  $G(3)/G(6)$ and $G(2)/G(6)$ are not subgroups of $K$. Hence $K$ is of geometric
   level 6.
  By (ii) of Theorem 4.1, $K$ has algebraic level 6
   (2 and 3 are primes in $\Bbb Z[\lambda]$). $K$ has 5 conjugates and one
   of them is given as follows. They are the only ones among the 26 groups with geometric
    level 6.
   $${\small
 -\infty _{_{_{\smile}}} \  \hspace{-.37cm} _{_{_{_{_{_{_{_{ \circ  }}}}}}}}
 0 \, _{_{_{\smile}}} \  \hspace{-.37cm} _{_{_{_{_{_{_{_{1 }}}}}}}}
\, 1/\lambda _{_{_{\smile}} }\   \hspace{-.37cm} _{_{_{_{_{_{_{_{\circ }}}}}}}}
  \lambda/\lambda _{_{_{\smile}} }\  \hspace{-.37cm} _{_{_{_{_{_{_{_{1 }}}}}}}}
  \lambda _{_{_{\smile}} }\  \hspace{-.37cm} _{_{_{_{_{_{_{_{\circ }}}}}}}}  \infty }\eqno(5.6)
 $$

\subsection{}Index 5,  $v_2 = 5$ with geometric level 5. Such subgroup is unique.
 A set of independent generators is given as follows.We will investigate
 this group further in a separate article (see [LL]).
    $${\small
 -\infty _{_{_{\smile}}} \  \hspace{-.37cm} _{_{_{_{_{_{_{_{\circ  }}}}}}}}
 0 \, _{_{_{\smile}}} \  \hspace{-.37cm} _{_{_{_{_{_{_{_{\circ }}}}}}}}
\, 1/\lambda _{_{_{\smile}} }\   \hspace{-.37cm} _{_{_{_{_{_{_{_{\circ}}}}}}}}
  \lambda/\lambda _{_{_{\smile}} }\  \hspace{-.37cm} _{_{_{_{_{_{_{_{\circ }}}}}}}}
  \lambda _{_{_{\smile}} }\  \hspace{-.37cm} _{_{_{_{_{_{_{_{\circ }}}}}}}}  \infty} \eqno(5.7)
 $$
%  If $K$  is congruence, by Proposition 5.1, $G(5) \subseteq K$.
% Since $G_5/G(2+\lambda ) \cong A_5$ is simple (see [LLT4]), $G_5/G(2+\lambda)$ does not have a normal
%  subgroup of index 5. It follows that $K$ is not a subgroup of $G(2+\lambda)$. This implies
%   that $KG(2+\lambda) = G_5$. As a consequence,
%   $A_5 \cong G_5/G(2+\lambda) \cong K/G(2+\lambda)\cap K$.
%    This implies that $G(2+\lambda)\cap K$ is a normal subgroup of $G_5$
%     of index $300 = 5|A_5|$. Note that
%    $$G(5) \triangleleft G(2+\lambda)\cap K \triangleleft G(2+\lambda)
%    \triangleleft G_5.\eqno(5.8)$$

\smallskip

\noindent {\bf Discussion.} It is well known that subgroups  of
 index $\le 6$ of the modular group $PSL(2, \Bbb Z)$
 are congruence. This is not true for the Hecke group $G_5$ as the
  above examples indicate (subsection 5.5).

 \section {Appendix A}
\noindent  Let $m\in \Bbb N$ be fixed. Set
 $${\tiny  A =T^m =  I + m\left ( \begin{array}{rc}
0 &\lambda \\
0& 0  \\
\end{array}
\right ),
 B = ST^mS^{-1}= I + m\left ( \begin{array}{rc}
0 &0 \\
-\lambda & 0  \\
\end{array}
\right ).}\eqno(A1)$$

 \noindent
Let $S$ and $T$ be given as in (1.1).
Denoted by
 $\Phi _m$
 the smallest normal subgroup of $G_5$ that contains $T^m$. Note that $T^{-m}\in \Phi_m$.
 The main purpose of this appendix is to find members of  $\Phi_m$.
  It is clear that
   $C= TST^mS^{-1}T^{-1} $ and $D = T^{-1}ST^{-m}S^{-1}T $ are members of $\Phi_m$.
 $C$ and $D$ are given as follows.
$${\tiny
C =
\left ( \begin{array}{rc}
1-m\lambda^2  & m\lambda^3 \\
-m\lambda& 1+m\lambda^2  \\
\end{array}
\right ),
D=
\left ( \begin{array}{rc}
1-m\lambda^2  & -m\lambda^3 \\
m\lambda& 1+m\lambda^2  \\
\end{array}
\right ) }
.\eqno(A2)$$
\noindent Note that $\lambda ^2 = \lambda +1$ and $\lambda^3= 2\lambda+1.$
 Hence the actual matrix form of $C$ and $D$ are given as follows as well.
 $${\tiny
  C= I + m\left ( \begin{array}{rc}
-\lambda -1  &2\lambda +1\\
-\lambda & \lambda +1  \\
\end{array}
\right ),
  D= I + m\left ( \begin{array}{rc}
-\lambda -1  &-2\lambda -1\\
\lambda & \lambda +1  \\
\end{array}
\right )}
.\eqno(A3)$$
Let  $p$ be an odd prime  divisor of $m$. Then  $ (I+mU)(I+mV) = I+ m(U+V)$ (mod $mp$). This transforms
 matrix multiplication  $(I+mU)(I+mV) $  modulo $mp$  into matrix addition $I+ m(U+V)$.
As a consequence, we have the following.
{\small $$  E = A^{-2}  B^{-1}C \equiv I + m {\tiny \left ( \begin{array}{rc}
-\lambda -1  &1\\
0 & \lambda +1  \\
\end{array}
\right ),}\,
F = A^2  BD  \equiv I + m   \left ( \begin{array}{rc}
-\lambda -1  &-1\\
0 & \lambda +1  \\
\end{array}
\right )
\,\,(\mbox{mod } mp)
\eqno(A4)$$}
\noindent
Note that $A, B, C, D, E,$ and $F$ are members of $\Phi_m$.
 We now conjugate $E$ and $F$ by $S$. It follows easily that
{\small
$$G =  S ES^{-1}  \equiv I + m {\tiny \left ( \begin{array}{rc}
\lambda +1  &0\\
-1 & -\lambda -1  \\
\end{array}
\right ),}\,
H =  S FS^{-1}  \equiv I + m {\tiny\left ( \begin{array}{rc}
\lambda +1  &0\\
1 & -\lambda -1  \\
\end{array}
\right )}\,\,(\mbox{mod } mp)
.\eqno(A5)$$}
\noindent Note that $G, H\in \Phi_m$.
Applying the identity
 $ (I+mU)(I+mV) \equiv  I+ m(U+V)$ (mod $mp$), we have
 $$GE \equiv
I + m {\tiny \left ( \begin{array}{rc}
 0 &1\\
-1 & 0 \\
\end{array}
\right ),}\,
GF \equiv
I + m {\tiny \left ( \begin{array}{rr}
 0 &-1\\
-1 & 0 \\
\end{array}
\right )}\,\,(\mbox{mod } mp)
\eqno(A6)$$
Recall that $p$ is odd.
Let $r \in \Bbb N$ be chosen such that $2r\equiv 1$ modulo $p$. It follows that
 $\Phi_m$ contains the following element.
$$ X = (GEGF)^r
\equiv
{\tiny \left (I + m \left ( \begin{array}{rc}
 0 &0\\
-2 & 0 \\
\end{array}
\right )\right ) ^r }
=
I + m {\tiny \left ( \begin{array}{rc}
 0 &0\\
 -1 & 0 \\
\end{array}
\right )}
={\tiny
\left ( \begin{array}{rc}
 1 &0\\
- m & 1 \\
\end{array}
\right )}
.
\,\,(\mbox{mod } mp)
\eqno(A7)$$
\noindent
{\bf Lemma A.} {\em Let $m\in\Bbb N$ and let $p$ be an odd prime divisor of $m$.
 Denoted by $\Phi_m$ the smallest normal subgroup of $G_5$ that contains $T^m$.
  Then $\Phi_m$ contains the set
$$\Delta_m = \{ T^m, ST^mS^{-1}, TST^mS^{-1}T^{-1},SXS^{-1}, X,
TSXS^{-1}T^{-1} \}.\eqno(A8)$$}

\noindent {\bf Remark.} The matrix representatives of the members in the set $\Delta_m$ modulo $mp$
are given as follows.
{\tiny $$ \left (\begin{array}{cc}
   1 & m\lambda\\
   0 & 1
 \end{array}\right ),
 \left (\begin{array}{cc}
   1 & 0\\
  - m\lambda & 1
 \end{array}\right ),
 \left  (\begin{array}{cc}
   1-m\lambda^2  & m\lambda^3\\
   -m\lambda& 1+m\lambda^2
 \end{array}\right ),
 \left (\begin{array}{cc}
   1 & m\\
   0 & 1
 \end{array}\right ),
 \left (\begin{array}{cc}
   1 & 0\\
  - m & 1
 \end{array}\right ),
 \left  (\begin{array}{cc}
   1-m \lambda & m\lambda^2 \\
  - m & 1+m\lambda
 \end{array}\right ).$$}
 Note that $\lambda^n= F_n \lambda+F_{n-1}$ where $F_k$ is the $k$-th Fibonacci numbers
  ($F_1=1, F_2=1, F_3=2, F_{k+1} =F_k+F_{k-1}$).
\section {Appendix B}
Let  $m$ be  an odd  rational integer
  and let $\Phi_{2m}$
  be the smallest normal subgroup of $G_5$ that contains $T^{2m}$.
   Following step by step  of what we have done in Appendix A ($(A1)$ to $(A6)$),
   $\Phi_{2m}$ consists of the following elements, $A$, $B$, $C$, $GE$ and $T(GE)T^{-1}$.
   The
     matrix representatives of these five matrices  modulo $4m$ are
{\tiny $$A_1=  I+ 2m\left (\begin{array}{cc}
   0 & \lambda\\
   0 & 0
 \end{array}\right ),A_2=  I+2m
 \left (\begin{array}{cc}
   0 & 0\\
  -\lambda & 0
 \end{array}\right ),A_3= I+2m
 \left  (\begin{array}{cc}
   -\lambda -1 &1\\
   -\lambda& \lambda +1
 \end{array}\right ), A_4= I+2m
 \left (\begin{array}{cc}
   0 & 1\\
   -1 & 0
 \end{array}\right )$$}
 and $ A_5= {\tiny I+2m
  \left (\begin{array}{cc}
   -\lambda  & \lambda \\
   -1 & \lambda
 \end{array}\right )}
 $. In order to determine the order of the group generated by these five matrices modulo $4m$,
  we consider $A_6 = A_1A_5$ and $A_7= A_1A_5A_2A_3$.
   The matrix representatives of $A_1, A_2, A_6, A_4$ and $A_7$ (in this order) are given as follows.
  {\tiny $$ I+ 2m\left (\begin{array}{cc}
   0 & \lambda\\
   0 & 0
 \end{array}\right ),  I+2m
 \left (\begin{array}{cc}
   0 & 0\\
  -\lambda & 0
 \end{array}\right ),
  I+2m
 \left (\begin{array}{cc}
   -\lambda  & 0\\
   -1 & \lambda
 \end{array}\right ),
 I+2m
 \left  (\begin{array}{cc}
   0 &1\\
   -1& 0
 \end{array}\right ),
  I+2m
  \left (\begin{array}{cc}
   -1  & 1 \\
   -1 & 1
 \end{array}\right )
 $$}
 Denoted by $E$ the group generated by the above five elements modulo $4m$.
 Applying the identity $(I+2mU)(I+2mV) \equiv I+2m(U+V)$ modulo $4m$, ones sees that
  the (1,1) and (2,2)-entries the matrix generated by the first three matrices
   is of the form $1+ 2m\lambda$ or 1 and
    the (1,1) and (2,2)-entries the matrix generated by the last two matrices
   is of the form $1+ 2m$ or 1. Consequently,
  one can show  that $E$ is elementary abelian of
  order $2^5$  if $m \ge 2$. In the case $m=1$, the last two matrices $A_4$ and $A_7$
   satisfy $A_4\equiv -A_7$ (mod 4). As we are working on $G_5$, where a matrix is
   identified with its negative, $E$ has order $2^4$.

\section {Appendix C}
Let $m$ be a rational integer. Suppose that $m$ is  a multiple of 4.
  One sees easily that
 $$(I+mU)(I+mV) \equiv I +m(U+V)\,\,\,  (\mbox{mod } 4m). \eqno(C1)$$
Let $S$ and $T$ be given as in (1.1) and construct matrices  $A, B, C, D, E, F$ as
 what we have done in Appendix A. It follows that
 $$
 GE \equiv
I + m \left ( \begin{array}{rc}
 0 &1\\
-1 & 0 \\
\end{array}
\right ),\,\,
GF \equiv
I + m \left ( \begin{array}{rr}
 0 &-1\\
-1 & 0 \\
\end{array}
\right )\,\,\, (\mbox{mod } 4m).\eqno(C2)$$
Hence
$$ {\small Y = GEGF
\equiv
I + m \left ( \begin{array}{rc}
 0 &0\\
-2 & 0 \\
\end{array}
\right )=
I + 2m \left ( \begin{array}{rc}
 0 &0\\
-1 & 0 \\
\end{array}
\right )
\,\,\, (\mbox{mod } 4m)}
.\eqno(C3)$$
Let $\Phi_{2m}$ be the smallest normal subgroup that contains  $T^m$.
 Similar to Appendix A, $\Phi_{2m}$ contains
$$
\Delta_{2m} = \{ T^{2m}, ST^{2m}S^{-1}, TST^{2m}S^{-1}T^{-1},SYS^{-1}, Y,
TSYS^{-1}T^{-1}\}. \eqno(C4)$$
 Similar to Appendix A,
 the actual matrix representatives of members in $\Phi_{2m}$  modulo $4m$ are given as follows.
{\tiny $$ \left (\begin{array}{cc}
   1 & 2m\lambda\\
   0 & 1
 \end{array}\right ),
 \left (\begin{array}{cc}
   1 & 0\\
  - 2m\lambda & 1
 \end{array}\right ),
 \left  (\begin{array}{cc}
   1-2m\lambda^2  &2 m\lambda^3\\
   -2m\lambda& 1+2m\lambda^2
 \end{array}\right ),
 \left (\begin{array}{cc}
   1 & 2m\\
   0 & 1
 \end{array}\right ),
 \left (\begin{array}{cc}
   1 & 0\\
  - 2m & 1
 \end{array}\right ),
 \left  (\begin{array}{cc}
   1-2m \lambda & 2m\lambda^2 \\
  - 2m & 1+2m\lambda
 \end{array}\right ).$$}
\section {Appendix D. Wohlfahrt's Theorem for Modular group}

The main purpose of this appendix is to give an elementary proof of Wohlfahrt's Theorem
 where Dirichlet's Theorem (see pp. 34 of [Wa]) is not used.

\smallskip
\noindent {\bf Lemma D1.} {\em  Let $p$ be a prime and let $A$ be a subgroup of $SL(2, p)$
 that contains at least two Sylow $p$-subgroups. Then $A = SL(2, p)$.}

\smallskip
\noindent {\em Proof.} Recall first that $|SL(2, p)| = (p+1)p(p-1)$ and that
 $SL(2, p)$ has $p+1$ Sylow $p$-subgroups.
 Since $A$ has at least two Sylow $p$-subgroups,
  Sylow's Theorem implies that  $A$ contains all the Sylow $p$-subgroups of $SL(2, p)$.
 In particular, $A$ contains all the elements of the following forms.
{\small  $$ U(x) =\left (
\begin{array}{rc}
1 & x \\
0 & 1 \\
\end{array}
\right ),
L(y)= \left (
\begin{array}{rc}
1 & 0 \\
y & 1 \\
\end{array}
\right ),\,\,\, x, y \in \Bbb Z_p.\eqno(D1)$$}
 Since $\Bbb Z_p$ is a field, one sees easily
 by applying elementary row and column operation that every element in $SL(2, p)$
  can be written as a word in $U(x)$ and $L(y)$. Hence $A= SL(2, p)$.\qed

\smallskip
\noindent
{\bf Lemma D2.} {\em
  Let $p, m \in \Bbb N $, where $p$ is a be a prime.
  Then $N(\Gamma (mp), T^m) = \Gamma (m)$,
  where  $N(\Gamma (mp), T^m) $ is the smallest normal subgroup of $SL(2, \Bbb Z)$
   that contains $\Gamma (mp)$ and $ T^m) $.
  }

    \smallskip
    \noindent {\em Proof.}
    We first assume that $p$ is a divisor of $m$.
     The index formula of $\Gamma (n)$ implies that  $|\Gamma(m)/\Gamma (mp)|= p^3$.
    Every element in $\Gamma (m)$ takes the form $I + m\sigma$, where $\sigma$ is a
     $2 \times 2$ matrix.
    Since $p$ is a prime divisor of $m$, one has $(I+m\sigma)(I+m\tau)
     \equiv I +m(\sigma + \tau)$ (mod $mp$).
      Hence  $\Gamma (m)/\Gamma (mp)$ is abelian and every non-identity element has order $p$.
       Equivalently,
     $\Gamma (m)/\Gamma (mp) \cong \Bbb Z_p\times \Bbb Z_p\times \Bbb Z_p$. The following is clear.
  $$ {\small  \Gamma (m)/\Gamma (mp) \cong\left <
U = \left (
\begin{array}{rc}
1 & m \\
0 & 1 \\
\end{array}
\right ),
V = \left (
\begin{array}{cr}
1 & 0 \\
-m & 1 \\
\end{array}
\right ),
W = \left (
\begin{array}{rr}
1-m & m \\
-m & 1+m \\
\end{array}
\right )
\right >  .}
\eqno(D2)$$
One sees easily that $U= T^m$, $V = S T^m S^{-1}$ and
 $W = TVT^{-1}$. Since these three elements are conjugates of $T^m$,
  we conclude that $\Gamma (m)/\Gamma (mp) \subseteq N(\Gamma (mp), T^m)/\Gamma (mp)$. Hence
   $\Gamma (m) = N(\Gamma (mp), T^m)$.
 We now assume that gcd$\,(p,m)=1$.
  The index formula of $\Gamma (n)$ implies that $\Gamma (m)/\Gamma (mp) \cong SL(2, \Bbb Z_p)$.
  It is clear that $\left < T^m\right >$ and
   $\left < ST^m S^{-1}\right >$ are two different Sylow $p$-subgroups of $N(\Gamma(mp), T^m)/
   \Gamma (mp)$. By Lemma D1, $N(\Gamma(mp), T^m) = \Gamma (m)$.\qed

\smallskip

\noindent {\bf Wohlfahrt's Lemma.} {\em
Let $r, s \in \Bbb N$ be given. Then the smallest normal subgroup of $SL(2, \Bbb Z)$
 that contains $\Gamma (rs)$ and $T^s$ is $\Gamma (s)$.}

\smallskip
\noindent {\em Proof.}
Let $k$ be the smallest positive integer such that $\Gamma (k) \subseteq N(\Gamma (rs), T^s)$.
 Note that $N(\Gamma (rs), T^s)\subseteq \Gamma (s)$.
 Hence
   $s|k$. Suppose that $s<k$. We may write $k$ into $k= mp$,
   where $p$ is a prime and $s|m$.
     By Lemma D2, $   \Gamma(m) = N(\Gamma (mp), T^m) =  N(\Gamma (k), T^m) \subseteq N(\Gamma (rs), T^s)$. This contradicts the minimality of $k$. Hence $N(\Gamma (rs), T^s) = \Gamma (s)$\qed

\smallskip\noindent {\bf Wolhfahrt's Theorem.} {\em Let $K \subseteq PSL(2, \Bbb Z)$ be a congruence
 subgroup of geometric level $r$
 and algebraic level $ s$. Then $r=s$.
}

\noindent {\em Proof.}
 Set $\overline \Gamma (r) = \Gamma(r)\left <\pm I\right >/\left <\pm I\right >$.
 Since the geometric level of $K$ is $r$,  $xT^{r}x^{-1} \in K$ for all $x$.
Let $ N(\overline \Gamma(rs), T^r) $ be the smallest normal subgroup of $PSL(2, \Bbb Z)$
 that contains $\overline \Gamma(rs)$ and $  T^r$.
  By Wohlfahrt's Lemma, $\overline \Gamma (r) =N(\overline \Gamma(rs), T^r) \subseteq K$. It follows that $\overline \Gamma (r)\subseteq \overline \Gamma (s)$
   as $s$ is the algebraic level. Hence $s$ is a divisor of $r$.
   On the other hand, since $\overline \Gamma (s) \subseteq K$,
    the geometric level $r$ is a divisor of $s$.
In summary, $r=s$.\qed

\bigskip
{\small
%\noindent CHENG LIEN LANG\\
\noindent Department of Mathematics, I-Shou  University, Kaohsiung, Taiwan,
Republic of China.

\noindent   \texttt{cllang@isu.edu.tw}

\smallskip
\noindent Singapore 669608, Republic of Singapore.

\noindent \texttt{lang2to46@gmail.com}}

\medskip

%\noindent congruence-1-6.tex

\end{document}